\documentclass{article}
\usepackage{amssymb}
\usepackage{amsthm,epsf,graphics,wrapfig,picinpar,fancyhdr,caption2,amsmath}
 \numberwithin{equation}{section}
\theoremstyle{definition}
\newtheorem{Theorem}{Theorem}

\newtheorem{Definition}{Definition}
\newtheorem{Remark}{Remark}
\newtheorem{Corollary}{Corollary}

\title{Hydrodynamic limit for a $d$-dimensional open symmetric exclusion process}
\date{}
\author{ Zhengye Zhou\footnote{Texas A$\&$M University, United States of America.  EMAIL: zyzhou@tamu.edu}}
 
\begin{document}

\maketitle
\begin{abstract}
In this paper we focus on the open symmetric exclusion process with parameter $m$ (open SEP($m/2$)), which allows $m$ particles each site and has an open boundary. We generalize the result about hydrodynamic limit for the open SEP$(m/2)$ originally raised in Theorem 4.12 of \cite{kuan2019stochastic}. We prove that the hydrodynamic limit of the density profile for a $d-$dimensional open SEP$(m/2)$ solves the $(d+1)-$dimensional heat equation with certain initial condition and boundary condition.
\end{abstract}

\section{Introduction}

Different types of interacting processes with open boundary condition have been studied in last few decades \cite{Derr1992}  \cite{Popkov_1999} \cite{Proeme_2010}, where the boundaries were seen as particle reservoirs and sinks. The systems conserve particle number away from boundaries and  exchange particles across their boundaries.
The symmetric exclusion process with an open boundary that allows up to $m$ particles each site  except for the boundary (open SEP($m/2$)) enjoys duality  \cite{GKRV09} \cite{kuan2019stochastic}. In \cite{kuan2019stochastic}, duality between open SEP($m/2$) with different types of boundaries was used to get the hydrodynamic limit for density profile and height function of a one-dimensional open SEP($m/2$) on $\mathbb{Z}$. In this paper, we start from the  duality result given in \cite{kuan2019stochastic} and  look into results about hydrodynamic limits in $d$ dimension. We use the method of Laplace transform to show that the hydrodynamic limits are solutions of certain partial differential equations. To the best of the author’s knowledge, this method has not occurred in proving the hydrodynamics of interacting processes.

We start by defining open SEP($m/2$) and two types of boundaries.

\begin{Definition}(Open SEP($m/2$)) 
Suppose $\mathcal{G}$ is a countable set, $\partial\mathcal{G}$ is a subset of $\mathcal{G}$, and $p$ is a symmetric stochastic matrix on $\mathcal{G}$.  The symmetric exclusion process with parameter $m\in\mathbb{N}$ and open boundary $\partial\mathcal{G}$ is a continuous time Markov process on particle configurations on $\mathcal{G}$. And $m$ is the maximum number of particles allowed for each site in the interior  $\mathcal{G}^\circ:=\mathcal{G}-\partial\mathcal{G}$,   $k_y\in \{0,1,\cdots,m\}$ is the occupation number at site $y\in\mathcal{G}^\circ$, $\alpha_x\in[0,1]$ is the boundary parameter for each $x\in\partial\mathcal{G}$. The jump rate for a particle from $x\in\mathcal{G}$ to $y\in\mathcal{G}$ is defined by 
\begin{equation}
\left\{\begin{array}{ll}
p(x,y)\alpha_x\frac{m-k_y}{m},  &  \text{ if } x\in\partial\mathcal{G} \text{ and } y\in\mathcal{G}^\circ.\\
 p(x,y)(1-\alpha_y)\frac{k_x}{m}, & \ \text{if $x\in\mathcal{G}^\circ$ and $y\in\partial\mathcal{G}$}.\\
 p(x,y)\frac{k_x}{m}\frac{m-k_y}{m}, &\ \text{if $x\in\mathcal{G}^\circ$ and $y\in\mathcal{G}^\circ$}.\\
 0, &\ \text{if $x\in\partial\mathcal{G}$ and $y\in\partial\mathcal{G}$}.
\end{array}\right.
\end{equation}
\end{Definition}
\begin{Definition}
(a). If $\alpha_x=0$ for all $x\in\partial\mathcal{G}$, the boundary for open SEP($m/2$) is a sink boundary. Each site $x\in\partial\mathcal{G}$ is a sink that absorbs particles from $\mathcal{G}^\circ$.\\\
(b). If $0<\alpha_x\le 1$  for all $x\in\partial\mathcal{G}$, the boundary for open SEP($m/2$) is a reservoir boundary. Each site $x\in\partial\mathcal{G}$ is a reservoir with infinitely many particles.
\end{Definition}
\begin{Remark}
In the case $\alpha_x=0$, the jump rate from site $x\in\partial\mathcal{G}$ is $0$. Thus the boundary site $x$ is absorbing, once a particle jumps there, it stays there forever. When $0<\alpha_x\le 1$, the jump rate from $x\in\partial\mathcal{G}$ is independent of the occupation number at $x$.
\end{Remark}
In order to lighten the notation, we omit the dependence on the dimension $d$ in definitions. For example, the hitting time $\tau_r$ and $\tau_{a,b}$, the functions $\rho_t,\phi$ and $\mathcal{N}$ are defined in $d$-dimension. Throughout this paper, $||\cdot||$ denotes the Euclidean norm in the corresponding dimension.  Let $B(0,r)$ be the open ball with radius $r$ centered at 0 in $\mathbb{R}^d $, and $B(0,r)^c$ be the complement of $B(0,r)$.  

Define the interior of $\mathcal{G}_L$ as $\mathcal{G}_L^\circ:=\mathbb{Z}^d\cap B(0,\sqrt{L})^c$,  the boundary of $\mathcal{G}_L$ as
$\partial \mathcal{G}_L:=\{z|z\in \mathbb{Z}^d,||z||\le\sqrt{L}, z \ \text{is adjacent to a vertex in }  \mathcal{G}_L^\circ\}$, and $\mathcal{G}_L:=\mathcal{G}_L^\circ\cup\partial\mathcal{G}_L$. Let $p\left(x,x\pm e_k \right)=\frac{1}{2d}$, where $\{e_k\}_{1\le k\le d}$ is the standard orthonormal basis for $\mathbb{R}^d$.

Our main result is about the hydrodynamic limit of an open SEP($m/2$).

\begin{Theorem}\label{th: 1}
Let $\mathfrak{s}_t$ evolve as an open SEP($m/2$)  on $\mathcal{G}_L\subset\mathbb{Z}^d$, with $\alpha_x=\alpha$ for all  $ x\in\partial\mathcal{G}_L$, where $0<\alpha\le 1$. Set $\mathfrak{s}_0(x)=0$ for all $x\in\mathcal{G}_L^\circ$.
Let $\rho_t(x)$ be the density profile of $\mathfrak{s}_t$, i.e.
\begin{equation}
    \rho_t(x)=\frac{1}{m}\left(\mathbb{P}(\mathfrak{s}_t(x)=1)+2\cdot \mathbb{P}(\mathfrak{s}_t(x)=2)+\cdots+m\cdot \mathbb{P}(\mathfrak{s}_t(x)=m)\right).
\end{equation}
Then the hydrodynamic limit $\phi(\chi,\tau):=\displaystyle\lim_{L\rightarrow\infty}\rho_{2dm\tau L}\left(\lfloor \chi L^{1/2}\rfloor\right)$ solves the heat equation
\begin{equation}\label{eq: 1a}
     \frac{\partial\phi(\chi,\tau)}{\partial \tau}=\frac{1}{2}\Delta \phi(\chi,\tau),
\end{equation}
on $\{\mathbb{R}^d-B(0,1)\}\times[0,\infty)$ with initial condition $\phi(\chi,0)=0$, and Dirichlet boundary condition $\phi(\chi,\tau)|_{||\chi||=1}=\alpha$. Here $\chi=(x_1,x_2\cdots,x_d)\in\mathbb{R}^d-B(0,1)$, and $\Delta=\displaystyle\sum_{i=1}^d \frac{\partial^2}{\partial x_i^2}$ is the $d-$dimensional Laplacian.
\end{Theorem}

\section{Proof of Theorem \ref{th: 1}}
First, let us recall the definition of stochastic duality and the duality result given by Theorem 4.11 in \cite{kuan2019stochastic}.
\begin{Definition}(Stochastic duality)
Two Markov processes $\mathfrak{s}_t$ and $\mathfrak{s}_t'$ on state spaces $\mathfrak{S}$ and $\mathfrak{S}'$ are dual with duality function $D$ on $\mathfrak{S}\times\mathfrak{S}'$ if
\begin{equation}
    E_{\mathfrak{s}}[D(\mathfrak{s}_t,\mathfrak{s}')]=E_{\mathfrak{s}'}[D(\mathfrak{s},\mathfrak{s}_t')] \text{ for all } \mathfrak{s} \in \mathfrak{S},\ \mathfrak{s}'\in \mathfrak{S}', \text{ and } t>0.
\end{equation}
On the left hand side, $E_{\mathfrak{s}}$ means $\mathfrak{s}_0=\mathfrak{s}$ and  on the right hand side, $E_{\mathfrak{s}'}$ means $\mathfrak{s}'_0=\mathfrak{s}'$.
\end{Definition}

\begin{Theorem}
Let $\mathfrak{s}_t$ evolve as an open SEP($m/2$) on $\mathcal{G}$ with a reservoir boundary $\partial\mathcal{G}$, and $\mathfrak{s}'_t$ evolve as an open SEP($m/2$) with a sink boundary $\partial\mathcal{G}$ and finitely many particles. Then $\mathfrak{s}_t$ and $\mathfrak{s}'_t$ are dual with respect to the function
\begin{equation}
   D(\mathfrak{s},\mathfrak{s}')= \prod_{y\in\partial\mathcal{G}} \alpha_y^{\mathfrak{s}'(y)}\prod_{x\in \mathcal{G}^\circ}\frac{\binom{\mathfrak{s}(x)}{\mathfrak{s}'(x)}}{\binom{m}{\mathfrak{s}'(x)}}1_{\{\mathfrak{s}(x)\ge\mathfrak{s}'(x)\}}.
\end{equation}
\end{Theorem}
Now we consider the special case of $\mathfrak{s}'_t$ that only consists of a single particle starting from site $x$. By observation, $\mathfrak{s}'_t$ is a simple random walk with jump rate $\frac{1}{2dm}$ that stops at the boundary $\partial\mathcal{G}_L$. Apply the duality relation, we have:
\begin{equation}
    \rho_t(x)=\mathbb{E}_\mathfrak{s}[D(\mathfrak{s}_t,\mathfrak{s}')]=\mathbb{E}_{\mathfrak{s}'}[D(\mathfrak{s},\mathfrak{s}'_t)]=\alpha\cdot\mathbb{P}_{x}(\inf_{0\le s\le t} ||\mathfrak{s}'_s||\le \sqrt{L}). \end{equation}
    Next, rescale time by $2dm$, so that the jump rate becomes $1$, we have
  \begin{equation}
  \rho_t(x)=\alpha\cdot\mathbb{P}_{x}(\inf_{0\le s\le t} ||\mathfrak{s}'_s||\le \sqrt{L}).  =\alpha\cdot\mathbb{P}_{x}(\inf_{0\le s\le \frac{ t}{2dm}} ||S_s||\le \sqrt{L}), 
\end{equation}
where $S_t$ is the $d$-dimensional continuous time random walks with jump rate 1, and initial condition $S_0=x$.

Now recall functional central limit theorem, which says $\frac{1}{\sqrt{L}}S_{tL}\Longrightarrow \mathcal{B}_t$, with $S_0=\lfloor\sqrt{L}\chi\rfloor$ and $\mathcal{B}_0=\chi$, where $\mathcal{B}_t$ is the standard $d-$dimensional Brownian motion. Thus, let  $t=2dm\tau L$, and $x=\lfloor \sqrt{L}\chi \rfloor$, with $||\chi||\ge 1$,
\begin{equation}\label{eq:2b}
\begin{aligned}
     &\phi(\chi,\tau)=\lim_{L\rightarrow \infty} \rho_{2dm\tau L}(\lfloor \sqrt{L}\chi \rfloor)= \lim_{L\rightarrow \infty} \alpha\mathbb{P}_{\lfloor \sqrt{L}\chi \rfloor}(\inf_{0\le s\le \frac{ t}{2dm}} ||S_s||\le \sqrt{L})\\ 
    & =  \alpha\mathbb{P}_{\chi}(\inf_{0\le s\le \tau} ||\mathcal{B}_s||\le 1)
    =\alpha\mathbb{P}_{\chi}(\tau_1\le \tau),
\end{aligned}
\end{equation}
where $\tau_r$ is the hitting time of $B(0,r) $ by $\mathcal{B}_t$. 

From equation \eqref{eq:2b}, the initial condition and boundary condition for $\phi$ follow easily by taking $\tau=0$ and $||\chi||=1$.

Next, we introduce Bessel process with index $v$. When $v=\frac{d-2}{2}\in \{\mathbb{N}^+/2-1\}$, Bessel process with index $v$ is identical in law with the Euclidean norm of the $d-$dimensional Brownian motion. Now define $\tau_{a,b}$ as the first hitting time to $b$ of the Bessel process starting at $a$. With this notation, we have
\begin{equation}\label{eq: 2d}
\phi(\chi,\tau)=\alpha\mathbb{P}_{\chi}(\tau_1\le \tau )=\alpha\mathbb{P}(\tau_{||\chi||,1}\le \tau ).
\end{equation}

Also, from \cite{2011arXiv1106.6132H}, we know the Laplace transformation of $\mathbb{P}(\tau_{||\chi||,1}\le \tau )$  with respect to $\tau$ is
\begin{equation}
    \mathcal{L} [\mathbb{P}(\tau_{||\chi||,1}\le \tau )](\lambda)=||\chi||^{-v}\frac{K_v(||\chi||\sqrt{2\lambda})}{  \lambda K_v(\sqrt{2\lambda})},
\end{equation}
where $K_v(z)$ is the second kind modified Bessel function of index $v$, which has integral form
\begin{equation}
    K_v(z)=\pi^{-1/2}\Gamma(v+1/2)(2z)^v\int_0^\infty \frac{cos(t)dt}{(t^2+z^2)^{v+1/2}}.
\end{equation}
It is a solution of the modified Bessel differential equation
\begin{equation}\label{eq: 2c}
z^2\frac{d^2f}{dz^2}+z\frac{df}{dz}-(z^2+v^2)f=0.    
\end{equation}
Then, apply Laplace transform to the function
\begin{equation}
    g(||\chi||,\tau):= \frac{\partial\phi(\chi,\tau)}{\partial \tau}-\frac{1}{2}\frac{\partial^2\phi(\chi,\tau)}{\partial ||\chi||^2}-\frac{d-1}{2||\chi||}\frac{\partial \phi(\chi,\tau) }{\partial ||\chi||}.
\end{equation}
Use the fact that for a differentiable function $f(t)$,
\begin{equation}
    \mathcal{L}[f'](s)=s\mathcal{L}[f](s)-f(0^-),
\end{equation}
and the modified Bessel differential equation \eqref{eq: 2c}, we get 
\begin{equation}\begin{aligned}
    &\mathcal{L}[g(||\chi||,\tau)](\lambda)=\frac{||\chi||^{-v-2}}{2\lambda K_v(\sqrt{2\lambda })}\left(-2\lambda||\chi||^2K''_v(\sqrt{2\lambda}||\chi||)\right.\\    & \left.-\sqrt{2\lambda}||\chi||K'_v(\sqrt{2\lambda}||\chi||)+(2\lambda ||\chi||^2+v^2)K_v(\sqrt{2\lambda}||\chi||)\right)=0.
\end{aligned}
\end{equation}
Thus, $g(||\chi||,\tau)=0$, i.e.
\begin{equation}
   \frac{\partial\phi(\chi,\tau)}{\partial \tau}=\frac{1}{2}\frac{\partial^2\mathcal{\phi}(\chi,\tau)}{\partial ||\chi||^2}+\frac{d-1}{2||\chi||}\frac{\partial \phi(\chi,\tau) }{\partial ||\chi||}. 
\end{equation}
Recall the $d-$dimensional Laplacian $\Delta$ in polar coordinates, write $\chi\in \mathbb{R}^d$  as $||\chi||\theta_{\chi}$,
\begin{equation}\label{eq: 2f}
    \Delta f(\chi) =\frac{\partial^2f(\chi)}{\partial ||\chi||^2}+\frac{d-1}{||\chi||}\frac{\partial f(\chi) }{\partial ||\chi||}+\frac{1}{||\chi||^2}\Delta_{S^{d-1}}f(\chi).
\end{equation}
Here $\Delta_{S^{d-1}}$ is the Laplace-Beltrami operator on the $(d-1)-$sphere. When $f$ is independent of $\theta$, $\Delta_{S^{d-1}}f=0$. From equation \eqref{eq: 2d}, it's  clear that function $\phi$ is radial in $\chi$, thus $\Delta_{S^{d-1}}\phi$   is zero.
Then, equation (\ref{eq: 1a}) follows.

%%%%%%%%%%%%%%%%%%%%%%%%%%%%%%%%%%%%%%%%%%%%%%
\section{Some applications}
In one-dimensional case, the hydrodynamic limit of the height function of the process $\mathfrak{s}_t$ is of great interest. The height function is usually defined as the number of particles to the right of a site on $\mathbb{Z}$ at a given time. In higher dimension, we can generalize it as the number of particles outside $B(0,r)$ at time $t$.
\begin{Corollary}\label{cor:1}
Define 
\begin{equation}
     N_r(\mathfrak{s}_t)=\sum_{||y||\ge r}\left(1_{\{\mathfrak{s}_t(y)=1\}}+2\cdot 1_{\{\mathfrak{s}_t(y)=2\}}\cdots+m\cdot 1_{\{\mathfrak{s}_t(y)=m\}} \right),
\end{equation}
 and
 \begin{equation}
  \mathcal{N}(r,\tau)=\lim_{L\rightarrow\infty}\mathbb{E}[m^{-1}N_{\sqrt{L}r}(\mathfrak{s}_{2dm\tau L})].
 \end{equation}
There exists a $M\in(0, \infty)$ such that $\mathcal{N}(r,\tau)$ solves the partial differential equation 
\begin{equation}\label{2}
    \frac{\partial\mathcal{N}(r,\tau)}{\partial \tau}=\frac{1}{2}\frac{\partial^2\mathcal{N}(r,\tau)}{\partial r^2}-\frac{d-1}{2r}\frac{\partial \mathcal{N}(r,\tau) }{\partial r}
\end{equation}
on $(1,\infty)\times[0,M]$ with initial condition $\mathcal{N}(r,0)=0$ and Neumann boundary condition $\frac{\partial \mathcal{N}(r,\tau)}{\partial r}|_{r=1}=-\frac{d\pi^{d/2}}{\Gamma(d/2+1)}\alpha$ when $d>1$, and   $\frac{\partial \mathcal{N}(r,\tau)}{\partial r}|_{r=1}=-\alpha$ when $d=1$.
\end{Corollary}

\subsection{Preliminaries }
Before proving Corollary \ref{cor:1}, we state some useful facts about functions $ \mathbb{P}(\tau_{r,1}\le \tau )$, $erfc(z)$ and $K_v(z)$ that will be used  in subsequent proof.
\begin{enumerate}
    \item 
We have uniform estimates for $\frac{\partial\mathbb{P}(\tau_{r,1}\le \tau )}{\partial \tau}$ when $r>1$ \cite{byczkowski2010hitting}:\\\ 
(a) For $d\ge 3$,
\begin{equation}\label{eq:3a}
    \frac{\partial\mathbb{P}(\tau_{r,1}\le \tau )}{\partial \tau}\approx \frac{r-1}{r}\frac{e^{-(r-1)^2/2\tau}}{\tau^{3/2}}\frac{1}{\tau^{(d-3)/2}+r^{(d-3)/2}}.
\end{equation}
(b) For $d=2$,
\begin{equation}\label{eq: 3b}
    \frac{\partial\mathbb{P}(\tau_{r,1}\le \tau )}{\partial \tau}\approx \frac{r-1}{r}e^{-(r-1)^2/2\tau}\frac{(r+\tau)^{1/2}}{\tau^{3/2}}\frac{1+\log r}{(1+\log (1+\frac{\tau}{r}))(1+\log (\tau+r))}.
\end{equation}
Here $f\approx g$ means that there exists strictly positive $c_1$ and $c_2$ depending only on $d$ such that $c_1 g\le f\le c_2 g$. 
\item
When $d=2$, there is another bound for $\mathbb{P}$ when $0<\tau<2r^2$ \cite{GRIGORYAN2002115}. There exists positive constant $ c_1, c_2$ such that 
\begin{equation}\label{eq: 3c}
   \mathbb{P}(\tau_{r,1}\le \tau )\le \frac{c_1}{\log r} e^{-
   \frac{c_2r^2}{\tau}}.
\end{equation}
\item
The complementary function $erfc(z)$ is defined by 
\begin{equation}
    erfc(z)=\frac{2}{\sqrt{\pi}}\int_z^\infty e^{-t^2}dt.
\end{equation}
When $||z||\rightarrow\infty$, it has asymptotic expansion  \cite{erfc}
\begin{equation}\label{eq: 3d}
\begin{aligned}
  erfc(z)=\frac{e^{-z^2}}{\sqrt{\pi}z}\left(1-\frac{1}{2z^2}+\cdots+(-1)^n\frac{(2n-1)!!}{(2z^2)^n}+\cdots\right)\\
  \thicksim1-\frac{\sqrt{z^2}}{z}+\frac{e^{-z^2}}{\sqrt{\pi}z}\left(1+O\left(\frac{1}{z^2}\right)\right).
\end{aligned}\end{equation}
\item
 When $|ph\ z|\le \frac{\pi}{2}$,  as $||z||\rightarrow\infty$,  $K_v(z)$ has asymptotic expansion \cite{Kv} 

\begin{equation}\label{eq: 3e}
    K_v(z)=\sqrt{\frac{\pi}{2z}}e^{-z}\left(\sum_{k=0}^\infty \frac{a_{k}(v)}{z^k}\right)=\sqrt{\frac{\pi}{2z}}e^{-z}\left(1+R_1(v,z)\right)
\end{equation}
with error bound $||R_1(v,z)||\le 2||\frac{4v^2-5}{z}||e^{\left|\left|\frac{v^2-\frac{1}{4}}{z}\right|\right|}$. 
\item
The derivative of $K_v(z)$ has the expression
\begin{equation}\label{eq: 3f}
    K'_v(z)=\frac{v}{z}K_v(z)-K_{v+1}(z).
\end{equation}
\end{enumerate}
\subsection{Proof of Corollary \ref{cor:1}}
Since the result for $d=1$ is already given by Theorem 4.12 in \cite{kuan2019stochastic}, we omit its proof. Next, we assume that the interchange of limits is justified in the following steps, some of which will be shown at the end of this section.

By definition of $\mathcal{N}$ and equation (\ref{eq:2b}), for $d>1$,
\begin{equation}
    \mathcal{N}(r,\tau)=\int_{||u||\ge r} \alpha\cdot \mathbb{P}(\tau_{||u||,1}\le \tau ) du=\frac{d\pi^{d/2}}{\Gamma(d/2+1)}\int_{r}^\infty w^{d-1}\alpha\cdot \mathbb{P}(\tau_{w,1}\le \tau ) dw.
\end{equation}
The second equation is obtained by changing variables to polar coordinates. Then the boundary condition and initial condition follow easily.

Next, we can write the partial derivatives of $\mathcal{N} $ in integral form as
\begin{equation}\label{eq: 3g}
    \frac{\partial\mathcal{N}(r,\tau)}{\partial \tau}
    =\frac{d\pi^{d/2}}{\Gamma(d/2+1)}\alpha\int_{r}^\infty \frac{\partial w^{d-1} \mathbb{P}(\tau_{w,1}\le \tau )}{\partial \tau}  dw;    
\end{equation}
\begin{equation}
    \frac{\partial \mathcal{N}(r,\tau) }{\partial r}=-\frac{d\pi^{d/2}}{\Gamma(d/2+1)}\alpha r^{d-1}\mathbb{P}(\tau_{r,1}\le \tau ) =\frac{d\pi^{d/2}}{\Gamma(d/2+1)}\alpha \int_{r}^\infty\frac{\partial w^{d-1} \mathbb{P}(\tau_{w,1}\le \tau )}{\partial w}  dw;
\end{equation}
 
\begin{equation}\label{eq: 3i}
    \frac{\partial^2\mathcal{N}(r,\tau)}{\partial r^2}=-\frac{d\pi^{d/2}}{\Gamma(d/2+1)}\alpha \frac{\partial r^{d-1}\mathbb{P}(\tau_{r,1}\le \tau )}{\partial r} =\frac{d\pi^{d/2}}{\Gamma(d/2+1)}\alpha \int_{r}^\infty\frac{\partial^2 w^{d-1} \mathbb{P}(\tau_{w,1}\le \tau )}{\partial w^2 } dw.
\end{equation}
\begin{Remark}
Equations \eqref{eq: 3g}-\eqref{eq: 3i} are derived  using the fundamental theorem of calculus, where some boundary terms are omitted in the above equations because they are zero, which will be shown later on.
\end{Remark}
Now, define function $\widetilde{g}(w,t)=w^{d-1} \mathbb{P}(\tau_{w,1}\le t )$. Then apply Laplace transformation to 
\begin{equation}
    H(w,t):= \frac{\partial \widetilde{g}(w,t)}{\partial t}-\frac{1}{2}\frac{\partial^2\widetilde{g}(w,t)}{\partial w^2}+\frac{d-1}{2w}\frac{\partial \widetilde{g}(w,t) }{\partial w}-\frac{d-1}{2w^2}\widetilde{g}(w,t), 
\end{equation}
we have 
\begin{equation}
\begin{aligned}
\mathcal{L}[H(w,t)](\lambda)=\frac{w^{v-1}}{2\lambda K_v(\sqrt{2\lambda})}\left((2w^2\lambda+v^2) K_v(w\sqrt{2\lambda})\right.\\
 \left.-2w^2\lambda K_v''(w\sqrt{2\lambda})-w\sqrt{2\lambda}K_v'(w\sqrt{2\lambda})\right).
\end{aligned}
\end{equation}
The modified Bessel equation \eqref{eq: 2c} again yields that 
\begin{equation}
    \mathcal{L}[H(w,t)]=0.
\end{equation}
Thus, $\widetilde{g}(w,t)$ satisfies the partial differential equation
\begin{equation}
  \frac{\partial \widetilde{g}(w,t)}{\partial t}=\frac{1}{2}\frac{\partial^2\widetilde{g}(w,t)}{\partial w^2}-\frac{d-1}{2w}\frac{\partial \widetilde{g}(w,t) }{\partial w}+\frac{d-1}{2w^2}\widetilde{g}(w,t).   
\end{equation}
Integrate with respect to $w$ on both sides from $r$ to infinity, the result for $\mathcal{N}$ follows.

To complete the proof, we need to show that the interchange of the integral and partial differential operators in equation \eqref{eq: 3g} is valid, and
\begin{equation}\label{eq:4a}
    \lim_{r\rightarrow\infty}r^{d-1}\mathbb{P}(\tau_{r,1}\le \tau )=0,
\end{equation}
\begin{equation}
    \lim_{r\rightarrow\infty}\frac{\partial r^{d-1}\mathbb{P}(\tau_{r,1}\le \tau )}{\partial r} =0.
\end{equation}

First, since $\mathbb{P}(\tau_{r,1}\le \tau )$ is monotone and continuous in $\tau$, we can pick $0<M<\infty$ such that $ \mathcal{N}(r,M)<\infty$. Then it suffices to show that for any fixed $r>1$, $\int_r^\infty\frac{\partial w^{d-1} \mathbb{P}(\tau_{w,1}\le \tau )}{\partial \tau} dw$ is uniformly convergent in $\tau$ on the region $[0,M]$. We aim to use dominated convergence theorem to show the uniform convergence. First, we need to find a dominated function for $\frac{\partial w^{d-1} \mathbb{P}(\tau_{w,1}\le \tau )}{\partial \tau}$.

Using equations \eqref{eq:3a} and \eqref{eq: 3b},  we obtain upper bounds listed below.  

There exists positive constants $C_d$ and $C_d'$ depending only on $d$ such that when $d\ge 3$, 
\begin{equation}\label{eq: 4b}
\begin{aligned}
&\frac{\partial\mathbb{P}(\tau_{w,1}\le \tau )}{\partial \tau}  \le C_d \frac{w-1}{w^{(d-1)/2}}\frac{e^{-(w-1)^2/2\tau}}{\tau^{3/2}}\\
&\le C_d' \left\{
\begin{array}{ll}
\frac{1}{w^{(d-1)/2}(w-1)^2},  &\  \text{when}\ w<\sqrt{3M}+1 \\
 \frac{w-1}{w^{(d-1)/2}}\frac{e^{-(w-1)^2/2M}}{M^{3/2}}, & \ \text{when}\ w\ge\sqrt{3M}+1
\end{array}
\right.;
\end{aligned}
\end{equation}
And when $d=2$,
\begin{equation}\label{eq: 4c}
\begin{aligned}
&\frac{\partial\mathbb{P}(\tau_{w,1}\le \tau )}{\partial \tau}  \le C_d \frac{(w-1)(w+M)^{\frac{1}{2}}}{w}\frac{e^{-(w-1)^2/2\tau}}{\tau^{3/2}}\\
 &\le C_d' 
    \left\{
\begin{array}{ll}
\frac{(w+M)^{\frac{1}{2}}}{w(w-1)^2},  &\  \text{when}\ w<\sqrt{3M}+1 \\
 (w+M)^{\frac{1}{2}}\frac{e^{-(w-1)^2/2M}}{M^{3/2}}, & \ \text{when}\ w\ge\sqrt{3M}+1
\end{array}
\right..
\end{aligned}
\end{equation}
Since the right most functions in inequalities \eqref{eq: 4b} and \eqref{eq: 4c} are integrable, we use them as dominate functions. Thus, for any fixed $r>1$,    $\int_r^\infty\frac{\partial w^{d-1} \mathbb{P}(\tau_{w,1}\le \tau )}{\partial \tau} dw$ is uniformly convergent.

Next, we proceed to prove equation \eqref{eq:4a}. When $d=2$, we can apply the bound in inequality \eqref{eq: 3c} directly.
When $d\ge 3$, we use inequality \eqref{eq: 4b} again,
\begin{equation}
\begin{aligned}
    \mathbb{P}(\tau_{r,1}\le \tau )\le \int_0^\tau C_d \frac{r-1}{r^{(d-1)/2}}\frac{e^{-(r-1)^2/2t}}{t^{3/2}} dt\le C_d \int_0^\tau (r-1)\frac{e^{-(r-1)^2/2t}}{t^{3/2}} dt\\
    =C'_d \cdot erfc(\frac{r-1}{\sqrt{2\tau}}). \end{aligned}
\end{equation}
Use the expansion of function $erfc(z)$ in equation \eqref{eq: 3d},
we have that for any finite $k\ge 0$ and fixed $\tau>0$,
\begin{equation}
    \lim_{r\rightarrow\infty}r^{k}\mathbb{P}(\tau_{r,1}\le \tau )=0.
\end{equation}

As for $    \frac{\partial\mathbb{P}(\tau_{r,1}\le \tau )}{\partial r}$, we begin by investigating its Laplace transform
\begin{equation}
     \mathcal{L} \left[\frac{\partial\mathbb{P}(\tau_{r,1}\le \tau )}{\partial r}\right](\lambda)=\frac{\sqrt{2}}{r^v}\frac{K'_v(r\sqrt{2\lambda})}{  \sqrt{\lambda} K_v(\sqrt{2\lambda})}-\frac{v}{r^{v+1}}\frac{K_v(r\sqrt{2\lambda})}{  \lambda K_v(\sqrt{2\lambda})}.
\end{equation}
Next, write $K'_v(z)$ as in equation \eqref{eq: 3f},
 by observation, it suffices to show
\begin{equation}\label{eq: 4d}
  \lim_{r\rightarrow\infty}  r^{v+1}\mathcal{L}^{-1}\left[\frac{K_{v+1}(r\sqrt{2\lambda})}{\sqrt{\lambda}K_v(\sqrt{2\lambda})}\right]=0.
\end{equation}
Since all the zeros of $K_v(z)$ have negative real part \cite{zeroBessel}, the inverse Laplace transform could be written as 
\begin{equation}
    \mathcal{L}^{-1}\left[\frac{K_{v+1}(r\sqrt{2\lambda})}{\sqrt{\lambda}K_v(\sqrt{2\lambda})}\right]=\frac{1}{2\pi i}\int_{1-i\infty}^{1+i\infty}e^{t\lambda}\frac{K_{v+1}(r\sqrt{2\lambda})}{\sqrt{\lambda}K_v(\sqrt{2\lambda})} d\lambda.
\end{equation}

For any fixed $\lambda\in 1+i\mathbb{R}$, using the asymptotic expansion \eqref{eq: 3d}, we have
\begin{equation}
  \lim_{r\rightarrow\infty}  r^{v+1}\frac{K_{v+1}(r\sqrt{2\lambda})}{\sqrt{\lambda}K_v(\sqrt{2\lambda})}=0.
\end{equation}

Last, we can bound $\left|\left|\frac{K_{v+1}(r\sqrt{2\lambda})}{\sqrt{\lambda}K_v(\sqrt{2\lambda})}\right|\right|$ by the module of $e^{-(r-1)\sqrt{2\lambda}}$, then apply dominated convergence theorem again, equation \eqref{eq: 4d} is proved by taking the limit inside the integral.

Left to show that $\left|\left|\frac{K_{v+1}(r\sqrt{2\lambda})}{\sqrt{\lambda}K_v(\sqrt{2\lambda})}\right|\right|$ is bounded by  $\left|\left|e^{-(r-1)\sqrt{2\lambda}}\right|\right|$. Let $\lambda=1+iy$, $y\in\mathbb{R}$ and define 
\begin{equation}
    f(y):=\left|\left|\frac{\frac{K_{v+1}(r\sqrt{2\lambda})}{\sqrt{\lambda}K_v(\sqrt{2\lambda})}}{e^{-(r-1)\sqrt{2\lambda}}}\right|\right|=r^{-\frac{1}{2}}\left|\left|\frac{1+R_1(v+1,r\sqrt{2\lambda})}{\lambda^{\frac{1}{2}}(1+R_1(v,\sqrt{2\lambda}))}\right|\right|,
\end{equation}
where in second equation, equation \eqref{eq: 3e} is used. Also, with the upper bound listed below equation \eqref{eq: 3e} for $||R_1(v,z)||$, we know that as $||z||\rightarrow \infty$, $||R_1(v,z)||\rightarrow 0$. So, for any fixed $r$, $f(y)$ is a continuous function from $\mathbb{R}$ to  $\mathbb{R}$ with 
\begin{equation}
    \lim_{y\rightarrow -\infty}f(y)=\lim_{y\rightarrow \infty}f(y)=0,
\end{equation}
Thus, for each fixed $r$, there exists a finite $M_r$ such that $0<f(y)\le M_r$. Moreover, this upper bound $M_r$ could be chosen such that it is decreasing in $r$.
The fact that  $\left|\left|e^{t\lambda}e^{-(r-1)\sqrt{2\lambda}}\right|\right|$
is integrable finishes the proof.

\end{document}